\documentclass{article}
\usepackage[margin=1in]{geometry}
\usepackage[utf8]{inputenc}
\usepackage{mathtools}

\usepackage{appendix}
\usepackage{setspace}
\usepackage{lineno}
\usepackage[english]{babel}
\usepackage{amsmath}
\usepackage[dvipsnames]{xcolor}
\usepackage{sidecap}
\usepackage{hyperref}
\usepackage{nicefrac}
\usepackage{fancyhdr}
\usepackage{amssymb}
\usepackage{amsthm}
\usepackage{bm}
\usepackage{natbib}
\usepackage{soul,xcolor}
\usepackage{graphicx}
\usepackage{cancel}

\usepackage{setspace} 
\usepackage{calc}
\usepackage{authblk}
\newtheorem{theorem}{Theorem}
\newtheorem{lemma}[theorem]{Lemma}
\newtheorem{prop}[theorem]{Proposition}
\newtheorem{corol}[theorem]{Corollary}

\graphicspath{ {./figures/} }
\DeclareGraphicsExtensions{.pdf,.png}
\renewcommand{\geq}{\geqslant}
\renewcommand{\leq}{\leqslant}

\newcommand{\rcp}[1]{\frac{1}{#1}} 
 
\newcommand{\cE}{{\mathcal E}}
\newcommand{\cG}{{\mathcal G}} 
\newcommand{\cU}{{\mathcal U}} 
 
\newcommand{\bt}{\begin{theo}} 
\newcommand{\et}{\end{theo}} 
\newcommand{\bl}{\begin{lemma}} 
\newcommand{\el}{\end{lemma}} 
\newcommand{\bp}{\begin{prop}} 
\newcommand{\ep}{\end{prop}} 
\newcommand{\bc}{\begin{corol}}
\newcommand{\ec}{\end{corol}}
\newcommand{\bdf}{\begin{df}} 
\newcommand{\edf}{\end{df}} 
\newcommand{\brem}{\begin{rem}} 
\newcommand{\erem}{\end{rem}} 
\newcommand{\bnrem}{\begin{nrem}} 
\newcommand{\enrem}{\end{nrem}} 
\newcommand{\bex}{\begin{ex}} 
\newcommand{\eex}{\end{ex}} 
\newcommand{\bcor}{\begin{cor}} 
\newcommand{\ecor}{\end{cor}} 
\newcommand{\bncor}{\begin{ncor}} 
\newcommand{\encor}{\end{ncor}} 
\newcommand{\bpf}{\begin{proof}} 
\newcommand{\epf}{\end{proof}}

\newcommand{\Du}{\operatorname{D}_u}

\newcommand{\pdiff}[3]{\frac{\partial^{#3}#1}{\partial #2^{#3}}}

\title{Enumerative combinatorics \\ of unlabeled and labeled time-consistent galled trees}
\author{Lily Agranat-Tamir\thanks{Department of Biology, Stanford University, Stanford, CA, USA.}, 
Michael Fuchs\thanks{Department of Mathematical Sciences, National Chengchi University, Taipei, Taiwan.}, 
Bernhard Gittenberger\thanks{Department of Discrete Mathematics and Geometry, Technische Universität Wien, Austria.}, 
Noah A.~Rosenberg$^*$}
\date{April 18, 2025}

\begin{document}
\maketitle

\begin{abstract}
\noindent In mathematical phylogenetics, the \emph{time-consistent galled trees} provide a simple class of rooted binary network structures that can be used to represent a variety of different biological phenomena. We study the enumerative combinatorics of unlabeled and labeled time-consistent galled trees. We present a new derivation via the symbolic method of the number of unlabeled time-consistent galled trees with a fixed number of leaves and a fixed number of galls. We also derive new generating functions and asymptotics for labeled time-consistent galled trees.
\end{abstract}
\medskip

\section{Introduction}
 
Rooted binary trees are fundamental combinatorial structures in mathematical phylogenetics. They are used for representing many aspects of biological descent that takes place in time---for example, the evolution of species from earlier species, the relationships among sequences of the same gene across different organisms, and the divergences of populations within a species. The root of a tree represents a common ancestor in the past and the leaves represent contemporaneous entities of the same type, such as species, genes, or populations. Nodes represent the divergence of an ancestor into two distinct descendants. 

In studying evolutionary descent, however, some biological processes---among them admixture, horizontal gene transfer, and hybridization---involve the merging of two entities. Representation of such processes requires a generalization from phylogenetic trees to phylogenetic networks. When viewing phylogenetic networks as structures unfolding in time, edges merge as well as diverge. 

Here, we provide new combinatorial analyses for a constrained set of (rooted, binary) phylogenetic networks, namely the time-consistent galled trees, or equivalently, the normal galled trees. We have previously studied the enumerative combinatorics of \emph{unlabeled} time-consistent galled trees~\citep{MathurAndRosenberg23, AgranatTamirEtAl24b, AgranatTamirEtAl24a}. We provide here new derivations of generating functions---which we derived previously using recursion in \cite{AgranatTamirEtAl24b}---by the \emph{symbolic method}, recalling the results of our earlier asymptotic analyses. We consider time-consistent galled trees with a fixed number of leaves and a fixed number of galls, and time-consistent galled trees with a fixed number of leaves and the number of galls unconstrained. Only limited results have been available for \emph{labeled} time-consistent galled trees~\citep{CardonaAndZhang20, FuchsAndGittenberger24}; we provide parallel combinatorial analyses of both unlabeled and labeled time-consistent galled trees. 

\section{Brief survey of past enumerative results}

\subsection{Network classes} 

Many classes of phylogenetic networks have been defined and their enumerative combinatorics studied. As background to our new results, we note existing enumerative results on certain classes of phylogenetic networks that, when constraints are sequentially applied, produce the time-consistent galled trees. For definitions of network classes, we rely primarily on the survey of \cite{KongEtAl22}, which described relationships among many network classes and explained their relevance to biology. Most past investigations consider leaf-labeled or vertex-labeled networks rather than unlabeled networks, with labels corresponding to specific entities such as species or genes. We report results on both leaf-labeled and unlabeled networks, and ``labeled'' networks henceforth refer to \emph{leaf}-labeled networks. 

We consider networks and trees that are rooted and binary. A \emph{rooted phylogenetic network} is a directed acyclic graph that has four properties. (i) A unique node, the \emph{root node}, has in-degree 0 and out-degree 2. (ii) Edges are directed away from the root node. (iii) \emph{Leaf nodes} have in-degree 1 and out-degree 0. (iv) Nonleaf, nonroot nodes have in-degree 2 and out-degree 1 (\emph{reticulation nodes} or \emph{hybrid nodes}) or in-degree 1 and out-degree 2 (\emph{tree nodes}). All structures considered are ``non-plane'': the left--right order in which the children of a node are indicated is not important. 

\begin{figure}[tbp]
    \centering
    \includegraphics[width=\textwidth]{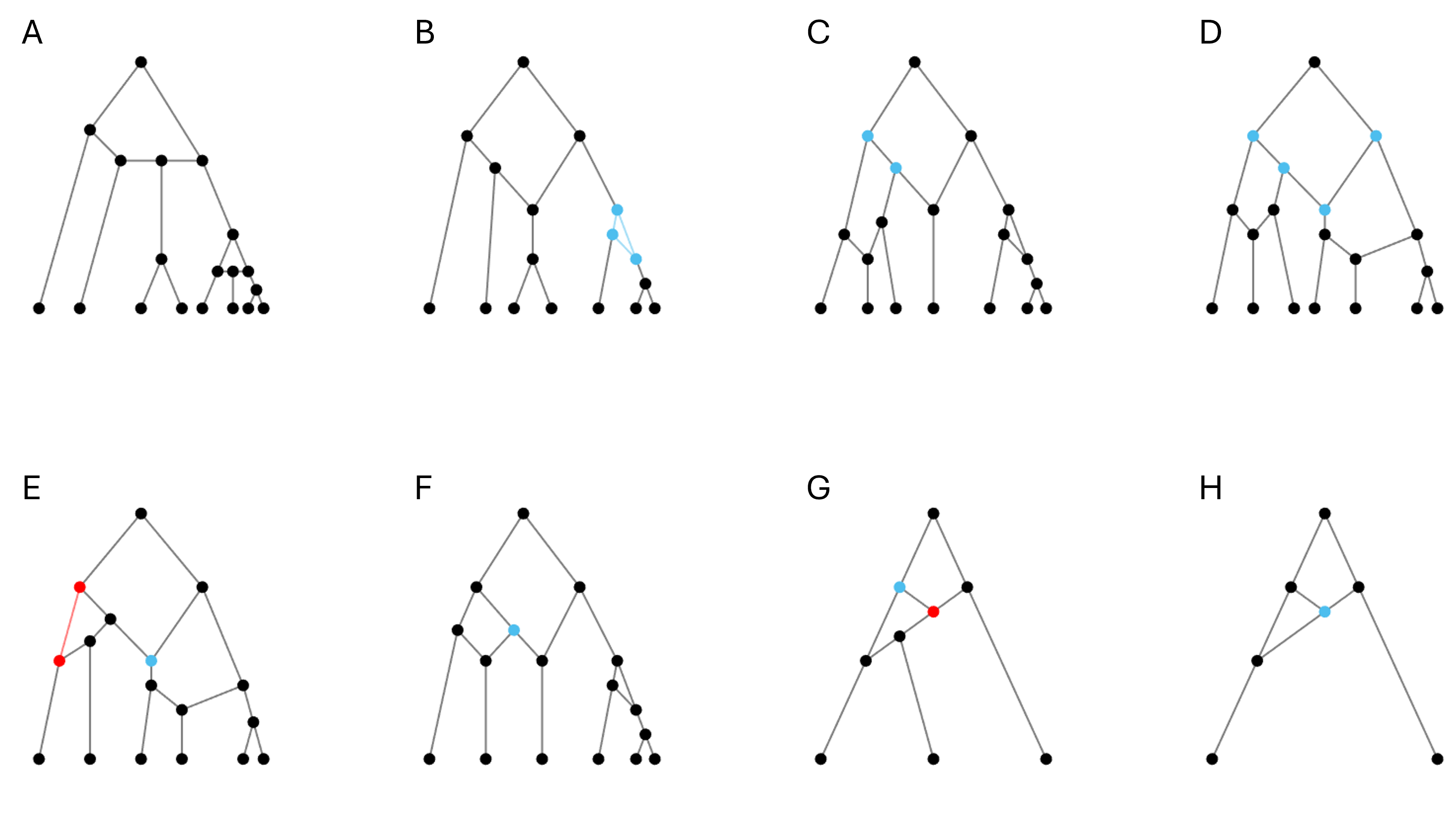}
    \vspace{-.7cm}
    \caption{Examples of phylogenetic networks. Sources for definitions appear in Table \ref{table:definitions}. 
    (A) Time-consistent galled tree (equivalently, normal galled tree). We draw reticulation events on a horizontal line to represent the concurrent existence of two merging entities that produce a hybrid entity. 
    (B) Galled tree. This network is not a time-consistent galled tree because it has a reticulation cycle in which the two parents of the reticulation are a parent--child pair themselves (blue). 
    (C) Galled tree-child network. This network is not a galled tree because it has nodes that are part of more than one reticulation cycle (blue). 
    (D) Normal network. This network is not a normal (time-consistent) \emph{galled tree} because it has nodes that are part of more than one reticulation cycle (blue).
    (E) Tree-child network. This network is not a galled tree-child network because it has a reticulation node that is in two reticulation cycles (e.g.~blue). It is not a normal network because it contains a ``shortcut'' (red).
    (F) Galled network. This network is not a galled tree-child network because it has a tree node that has only reticulation nodes as children (blue). 
    (G) Reticulation-visible network. This network is not a galled network because it has a reticulation node that is in two reticulation cycles (red). It is not a tree-child network because it has a tree node that is a parent of two reticulation nodes (blue). 
    (H) Phylogenetic network. This network is not a reticulation-visible network because it has a reticulation node all of whose descendant leaves possess paths from the root that do not traverse it (blue).} 
    \label{fig:1}
\end{figure}

\begin{figure}[tbp]
    \centering
    \vspace{-.5cm}
    \includegraphics[width=6.5cm]{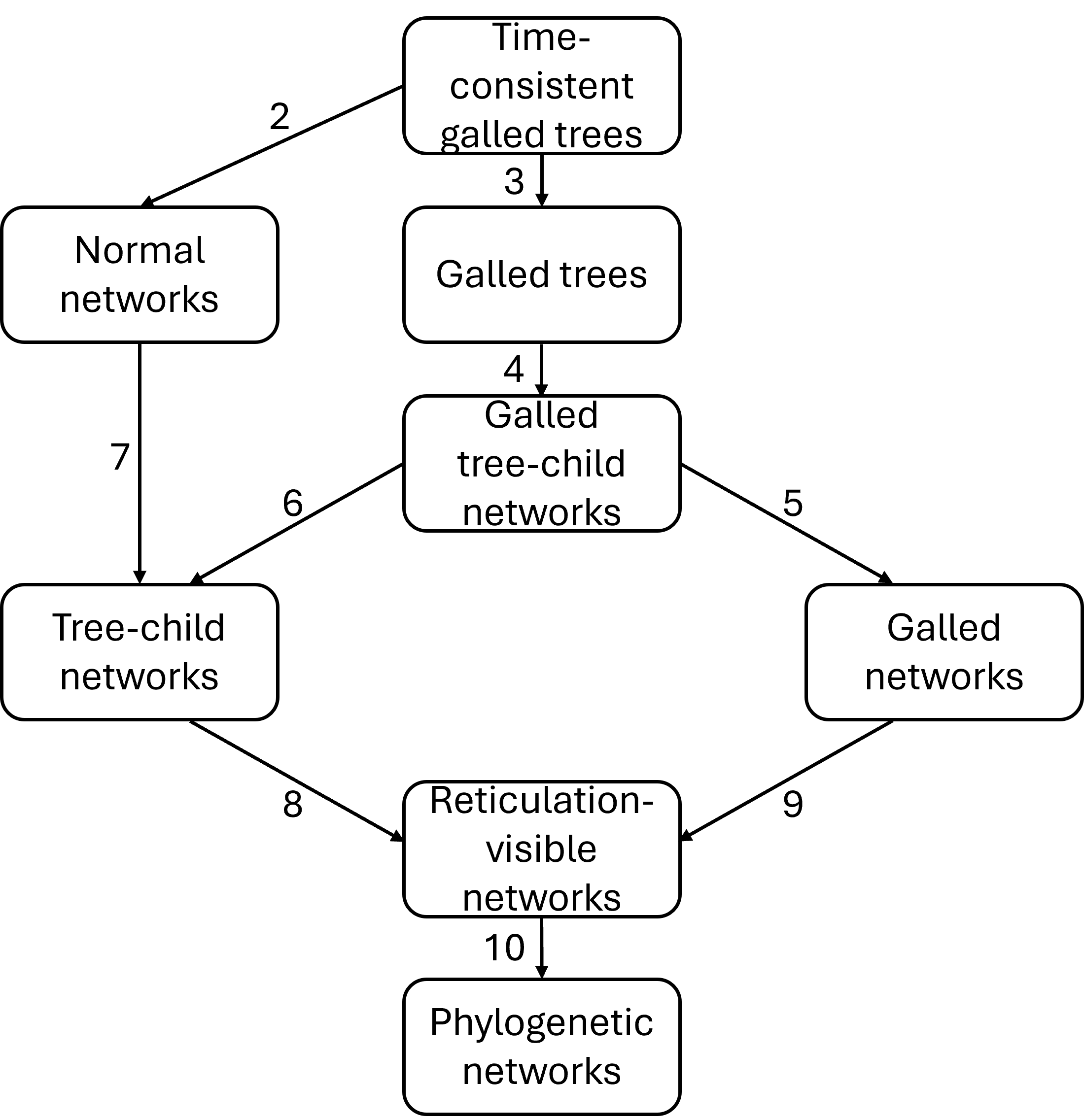}
    \vspace{-.35cm}
    \caption{Inclusion relations between classes of phylogenetic networks. Arrows represent the inclusion of the class on the top in the class below it. Notice that time-consistent galled trees are also normal galled trees because time-consistency implies no ``shortcuts'' (in a galled tree, a shortcut can only appear in a gall, contradicting time-consistency). Inclusion relationships, indicated by numbers, are described in Table \ref{table:inclusion}.} 
    \label{fig:2}
\end{figure}

We review enumerative results for several classes of networks: time-consistent galled trees (equivalent to normal galled trees, as we discuss), galled trees, galled tree-child networks, normal networks, tree-child networks, galled networks, reticulation-visible networks, and phylogenetic networks in general. We show illustrative examples of these network classes in Figure~\ref{fig:1}; sources for definitions appear in Table \ref{table:definitions}. The inclusion relations of the classes appear in Figure~\ref{fig:2}, and they are summarized in Table~\ref{table:inclusion}. 

The particular classes of networks that we survey are chosen for two reasons. First, the time-consistent galled trees---on which our new analysis focuses---represent a subset of each class, so that the various classes are meaningfully connected to the structures of primary interest. Second, we seek to compare results that we report on the number of time-consistent galled trees with a fixed number of galls to corresponding results on other classes of phylogenetic networks with a fixed number of reticulations; suitable results are available for the classes that we consider. 

Enumerative combinatorics results for classes of phylogenetic networks typically fall into one of the following categories: (1) exact counts of the number of networks with $n$ leaves, generating functions for these counts, and associated asymptotic approximations; (2) exact counts of the number of networks with $n$ leaves and a specified value $k$ for an additional parameter describing the number of reticulations, multivariate generating functions involving the number of leaves and the additional parameter, and associated asymptotic approximations for counts with the additional parameter fixed; and (3) detailed analyses of network enumerations with $n$ leaves and $k$ reticulations for special cases of $k$ (e.g.~$k=1,2$). References for results in these categories appear in Table \ref{table:survey}.

\clearpage

\clearpage 

\begin{table}[thb] 
\centering
\caption{Sources for definitions of network classes in Figure \ref{fig:1}.} 
\label{table:definitions} \scriptsize
\vspace{0.0cm}
\begin{tabular}{|l|l|l||}
\hline 
Network class                 & Source for definition \\ \hline
Time-consistent galled trees  & \cite{AgranatTamirEtAl24a}, p.~3\\
Galled trees                  & \cite{KongEtAl22}, p.~16 \\ 
Galled tree-child networks    & \cite{ChangEtAl24}, Definition 5 \\
Normal networks               & \cite{KongEtAl22}, p.~10\\
Tree-child networks           & \cite{KongEtAl22}, p.~9\\
Galled networks               & \cite{KongEtAl22}, p.~16 \\
Reticulation-visible networks & \cite{KongEtAl22}, p.~11\\
Phylogenetic networks         & \cite{KongEtAl22}, p.~5 \\
\hline
\end{tabular}
\end{table}

\begin{table}[thb] 
\centering
\caption{Inclusion relationships for categories of networks.} 
\label{table:inclusion} \scriptsize
\vspace{0.0cm}
\begin{tabular}{|c|c|c|c|c|}
\hline 
Line & First network class & Relationship & Second network class & Source \\ \hline
1 & Time-consistent galled trees & identical to & Normal galled trees & \cite{AgranatTamirEtAl24a}, p.~3 \\
2 & Time-consistent galled trees & proper subset of & Normal networks & \cite{AgranatTamirEtAl24a}, p.~3 \\
3 & Time-consistent galled trees & proper subset of & Galled trees & \cite{AgranatTamirEtAl24a}, p.~3 \\
4 & Galled trees                 & proper subset of & Galled tree-child networks & \cite{KongEtAl22}, Table 1 \#13 \\
5 & Galled tree-child networks   & proper subset of & Galled networks & \cite{ChangEtAl24}, Remark 4\\ 
6 & Galled tree-child networks   & proper subset of & Tree-child networks & \cite{ChangEtAl24}, Remark 4 \\
7 & Normal networks              & proper subset of & Tree-child networks & \cite{KongEtAl22}, Fig.~12\\
8 & Tree-child networks          & proper subset of & Reticulation-visible networks & \cite{KongEtAl22}, Fig.~12 \\
9 & Galled networks              & proper subset of & Reticulation-visible networks & \cite{KongEtAl22}, Fig.~12 \\
10 & Reticulation-visible networks &  proper subset of & Phylogenetic networks & \cite{KongEtAl22}, Fig.~5\\
\hline
\end{tabular}
\end{table}

\begin{table}[tbh]
\caption{Summary of enumerative results for phylogenetic networks}
\label{table:survey}\scriptsize
\begin{tabular}{|p{2.2cm}|p{2.9cm}|p{3cm}|p{3cm}|p{3.3cm}|} \hline
Network class & \multicolumn{4}{|c|}{Type of result} \\[1ex] \cline{2-5}
              & & \multicolumn{3}{|c|}{Leaf-labeled networks} \\[1ex] \cline{3-5}
              & Unlabeled networks & Fixed number of leaves $n$ & Fixed number of leaves $n$ and fixed number of reticulations $k$ & Fixed number of leaves $n$ and special cases of $k$ (e.g.~$k=1,2$) \\[1ex] \hline
Time-consistent galled trees & \cite{MathurAndRosenberg23}, \cite{AgranatTamirEtAl24a}, \cite{AgranatTamirEtAl24b}, this study &\cite{CardonaAndZhang20}, \cite{FuchsAndGittenberger24}, this study & This study & This study\\[1ex] \hline
Galled trees & & \cite{BouvelEtAl20}, \cite{CardonaAndZhang20} & \cite{BouvelEtAl20}, \cite{CardonaAndZhang20} & \cite{CardonaAndZhang20} \\[1ex] \hline
Galled tree-child networks & & \cite{ChangEtAl24} & \cite{ChangEtAl24} & \\[1ex] \hline
Normal networks & & &\cite{FuchsEtAl19}, \cite{FuchsEtAl21}, \cite{FuchsEtAl22b}, \cite{FuchsEtAl24} & \cite{CardonaAndZhang20}, \cite{FuchsEtAl21} \\[1ex] \hline
Tree-child networks & & \cite{FuchsEtAl21b}, \cite{BienvenuEtAl21} & \cite{FuchsEtAl19}, \cite{CardonaAndZhang20}, 
\cite{FuchsEtAl21},
\cite{BienvenuEtAl21},
\cite{FuchsEtAl22b} & \cite{CardonaAndZhang20}, \cite{FuchsEtAl21} \\[1ex] \hline
Galled networks & & \cite{GUNAWAN2020644}, \cite{FuchsEtAl22a}, \cite{ChangAndFuchs24} & \cite{ChangAndFuchs24} & \cite{CardonaAndZhang20}, \cite{ChangAndFuchs24} \\[1ex] \hline
Reticulation-visible networks & & \cite{ChangAndFuchs24} & \cite{ChangAndFuchs24} & \cite{CardonaAndZhang20}, \cite{ChangAndFuchs24} \\[1ex] \hline
Phylogenetic networks & & & \cite{Mansouri22} & \cite{CardonaAndZhang20}, \cite{Mansouri22} \\ \hline 
\end{tabular}
\end{table}

\subsection{Asymptotic equivalence of network classes with fixed numbers of reticulations}
\label{sec:review}

Several studies have examined numbers of phylogenetic networks in various classes, considering $n$ leaves and $k$ reticulations, for fixed $k$, as $n \rightarrow \infty$. Although a general phylogenetic network need not be reticulation-visible (Figure \ref{fig:1}), \cite{ChangAndFuchs24} have shown in their Corollary 1 that the number of labeled reticulation-visible networks $RV_{n,k}$ is asymptotically equivalent to the number of general phylogenetic networks $PN_{n,k}$, as examined by \cite{Mansouri22}. \cite{ChangAndFuchs24} in their Theorem 2 find that the number of galled networks $GN_{n,k}$ has the same asymptotic equivalence. \cite{FuchsEtAl22b} in their Theorem 1 show further that the number of tree-child networks $TC_{n,k}$ also has this asymptotic equivalence, and in their Corollary 2 that the number of normal networks $N_{n,k}$ does as well, a result explored in further detail in Corollary 7 of \cite{FuchsEtAl24}. Theorem 24 of \cite{ChangEtAl24} finds the same asymptotic equivalence for galled tree-child networks, $GTC_{n,k}$. In summary,
\begin{equation}
\label{eq:asymptotic}
    PN_{n,k} \sim RV_{n,k} \sim GN_{n,k} \sim TC_{n,k} \sim N_{n,k} \sim GTC_{n,k} \sim \frac{2^{k-1} \sqrt{2}}{k!} \bigg(\frac{2}{e}\bigg)^n n^{n+2k-1}.
\end{equation}

\subsection{Time-consistent galled trees}

The study of time-consistent galled trees begins with \cite{MathurAndRosenberg23}, who argued that these networks have a natural biological relevance, as they represent a relatively simple class of networks (galled trees) that can be understood via evolutionary processes unfolding in time (time-consistency). \cite{MathurAndRosenberg23} described how to count the number of labeled histories that are compatible with a time-consistent galled tree. \cite{MathurAndRosenberg23} algorithmically enumerated unlabeled time-consistent galled trees in a procedure that was formalized by \cite{AgranatTamirEtAl24a}. 

As the focus of \cite{MathurAndRosenberg23} had been on labeled histories for time-consistent galled trees, a subsequent study \citep{AgranatTamirEtAl24a} considered the enumerative combinatorics of the (unlabeled) time-consistent galled trees themselves. In that study, we obtained (1) a recursion for the exact count of the number of unlabeled time-consistent galled trees on $n$ leaves~\cite[eqs.~15 and 16]{AgranatTamirEtAl24a}; (2) a generating function for these counts~\cite[eq.~36]{AgranatTamirEtAl24a}; and (3) associated asymptotic approximations~\cite[eq.~42]{AgranatTamirEtAl24a}. We also obtained (4) a recursion for the exact count of the number of networks with $n$ leaves and $g$ galls~\cite[eqs.~26 and 27]{AgranatTamirEtAl24a}; (5) a bivariate generating function involving the number of leaves and the number of galls~~\cite[eq.~56]{AgranatTamirEtAl24a}. Subsequently, we studied (6) the asymptotic approximation to the number of time-consistent galled trees with a fixed number of galls $g$~\citep[Theorem 10]{AgranatTamirEtAl24b}. 

Left unanalyzed are the labeled time-consistent galled trees, for comparison to the enumerations of other classes of labeled phylogenetic networks discussed in Section \ref{sec:review}. Further, many of the results of \cite{AgranatTamirEtAl24a} and \cite{AgranatTamirEtAl24b} can be obtained efficiently using the symbolic method of \cite{FlajoletAndSedgewick09}. We proceed to these analyses. 

\section{Symbolic method for unlabeled time-consistent galled trees}
\label{sec:unlabeled}

\subsection{Overview}

Generating functions enumerating unlabeled time-consistent galled trees were obtained by recurrences in sections 5.1 and 5.3 of \cite{AgranatTamirEtAl24a} and sections 5.1 and 6.1 in \cite{AgranatTamirEtAl24b}; Section \ref{sec:unlabeled} here derives them by the symbolic method. In Section \ref{sec:labeled}, we consider \emph{leaf-labeled} time-consistent binary galled trees, deriving the exponential generating functions and asymptotic approximations of these networks with exactly one gall, with exactly two galls and with any fixed number of galls.

\subsection{Definitions}
\label{sec:unlabeled-definitions}

\emph{Time-consistent galled trees} are rooted binary phylogenetic networks with the following properties: (i) each reticulation node $a_r$ has a unique ancestor node $r$ with exactly two non-overlapping paths of edges connecting $r$ to $a_r$. Ignoring the direction of the edges, the two paths from $r$ to $a_r$ produce a cycle $C_r$: a \emph{gall}. (ii) For reticulation nodes $a_r$ and $a_s$, $a_r \neq a_s$, associated with galls $C_r$ and $C_s$, the sets of nodes in the associated galls are disjoint. (iii) Ancestor node $r$ and reticulation node $a_r$ are separated by two or more edges. This last condition encodes the requirement that we consider only \emph{normal} galled trees, equivalent to \emph{time-consistent} galled trees. We refer to the galled trees rooted at internal nodes of a galled tree as \emph{subtrees}. 

For the time-consistent galled trees, viewing galls as representations of biological merging events, we depict hybridizing nodes and their associated hybrid node on a horizontal line, representing the simultaneity of these nodes when a galled tree is taken to represent a structure evolving in time. This representation is possible in time-consistent networks and not in general galled trees, where the two parents of a reticulation node can have a parent--child relation as well. The definition of time-consistent galled trees follows \cite{MathurAndRosenberg23}, \cite{AgranatTamirEtAl24a}, and \cite{AgranatTamirEtAl24b}.

\subsection{Symbolic method}

We obtain generating functions via the \emph{symbolic method}. In this powerful approach to combinatorial enumeration, a class of structures is described in terms of a canonical series of constructions. A generating function enumerating objects in the class is obtained by use of a library that translates elements in the construction into algebraic expressions in the generating function. \cite{FuchsAndGittenberger24} summarized key components of the approach; a full description appears in Part A of \cite{FlajoletAndSedgewick09}.

\subsection{Generating functions}

\subsubsection{No galls}

To study generating functions for galled trees, first denote by $\cU(t)$ the generating function counting unlabeled trees with no galls. An unlabeled tree with no galls $\cU$ can be either a single leaf, $\{ \square \}$, or it can consist of a root $\{ \circ \}$ to which a multiset of two (possibly identical) unlabeled trees are attached, $\textrm{MSET}_2(\cU)$. According to the symbolic method, we have the construction
\begin{equation}
   \cU = \{ \square \} \textrm{   }\, {\dot{\cup}}  \,\textrm{   }
   \{ \circ \} \times
   \textrm{MSET}_2(\cU),
\end{equation}
where $\dot{\cup}$ represents a disjoint union (combinatorial sum, \cite[section 2.1]{FuchsAndGittenberger24}), and $\times$ represents a Cartesian product. The library \citep[Figure I.18, p.~93]{FlajoletAndSedgewick09} translates this construction into \citep[I.44, p.~72]{FlajoletAndSedgewick09} 
\begin{align}
\cU(t) = t + \frac{1}{2} [\cU(t)^2 + \cU(t^2)].
\end{align}

The asymptotic growth of the coefficient $U_n$ of $\cU(t)$, representing the number of unlabeled trees with $n$ leaves is well-known \citep[][p.~55]{Otter48, Harding71}. The form for $U_n$ in eq.~2 of \cite{AgranatTamirEtAl24b} is convenient:
\begin{equation}
U_n \sim \frac{\gamma}{2 \sqrt{\pi} }n^{-3/2} \rho^{-n},
\end{equation}
for constants $\gamma \approx 1.13000$, $\rho \approx 0.4027$, and $\sqrt{\pi}=\Gamma(\frac{1}{2})$. 

\begin{figure}[tb]
    \centering
    \includegraphics[width=4.5cm]{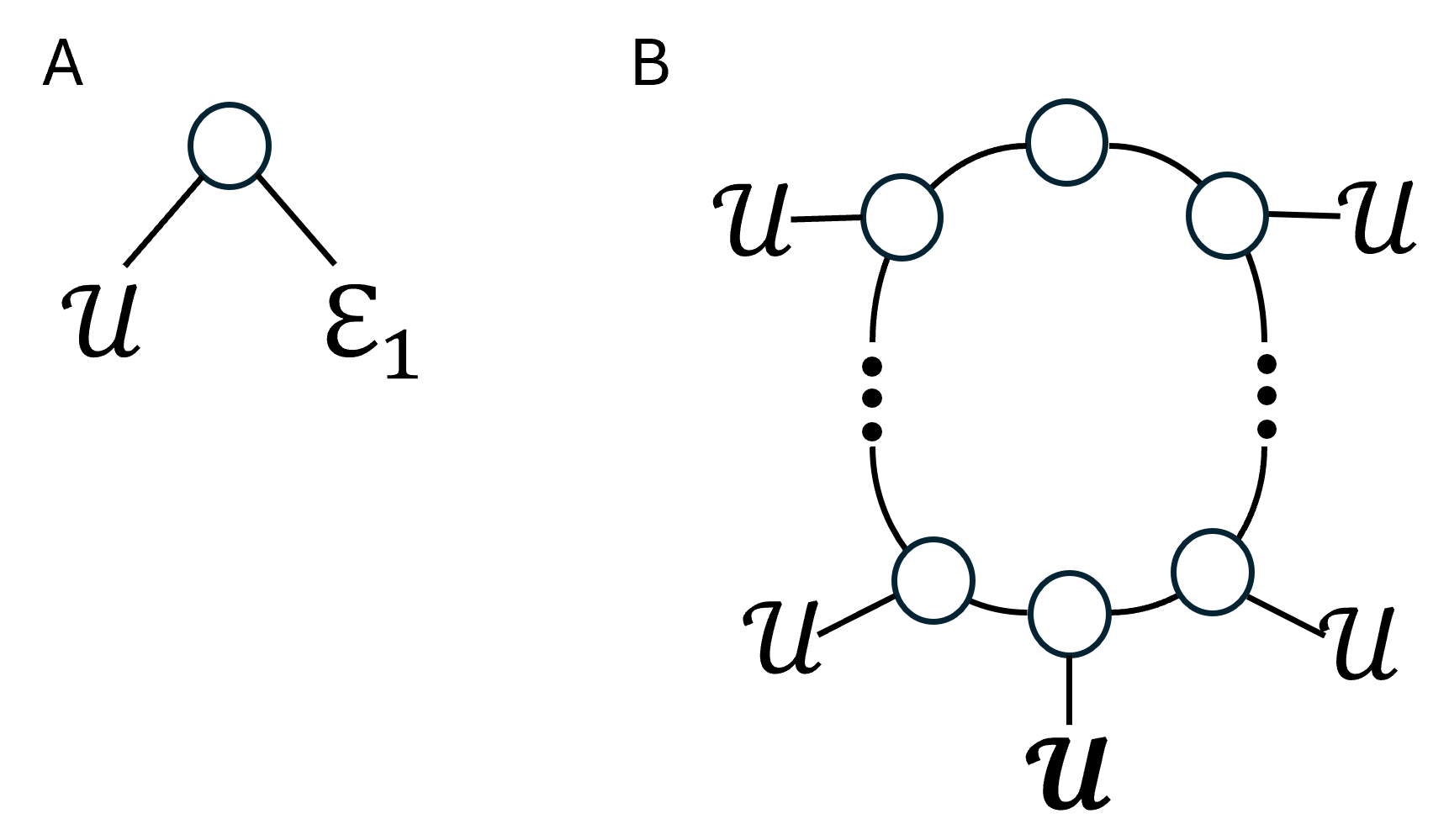}
    \caption{A non-plane unlabeled time-consistent galled tree with one gall, $\mathcal{E}_1$, has one of two structures. (A) A root with one subtree with no galls ($\mathcal{U}$) and one subtree with one gall ($\mathcal{E}_1$); there is no case of symmetry in this scenario. (B) A root gall, which is the only gall, and a subtree with no galls descended from the reticulation node (bold), with two non-empty sequences of subtrees with no galls (non-empty because of the time-consistency condition). The latter scenario has a case of symmetry, and because the trees are unlabeled in addition to being non-plane, the two non-empty sequences form a multiset of size 2.}
    \label{fig:E1}
\end{figure}

\subsubsection{One gall}

Next, we derive $\cE_1(t)$, the generating function counting the time-consistent galled trees with exactly one gall. If a galled tree has exactly one gall, then (1) it does not have a root gall and one of the root's subtrees contains the gall (Figure \ref{fig:E1}A), or (2) it has a root gall (Figure \ref{fig:E1}B). 

Using the terminology of the symbolic method, in the former case, the root node is attached to a galled tree with no galls $(\cU)$ and to a tree with one gall $(\cE_1)$. In the latter case, the structure contains a root node attached to an unordered pair $(\textrm{MSET}_2)$ of paths of unlabeled trees, with at least one node per path $(\textrm{SEQ}^+)$, together with the child of the reticulation node---an unlabeled tree as well. The construction gives
\begin{align}
    \mathcal{E}_1 =  \{ \circ \} \times  
    \Big[
    \underbrace{ \mathcal{U}  \times \mathcal{E}_1}_{(1)} 
    \textrm{   }\, \dot{\cup} \,\textrm{   }
    \underbrace{\mathcal{U} \times \textrm{MSET}_2   \Big( \textrm{SEQ}^+(\mathcal{U})\Big)}_{(2)}
    \Big].
\end{align}
Converting to a generating function~\citep[Figure I.18, p.~93]{FlajoletAndSedgewick09}, we have
\begin{align*}
\cE_1(t) &= \underbrace{ \cU(t) \, \cE_1(t)}_{(1)} + \underbrace{ \frac{\cU(t)}2\Bigg[\bigg(\frac{\cU(t)}{1-\cU(t)}\bigg)^2+\frac{\cU(t^2)}{1-\cU(t^2)}\Bigg]}_{(2)}.
\end{align*}
This expression simplifies to
\begin{align}
\cE_1(t) &= \frac{\cU(t)^3}{2[1-\cU(t)]^3}+\frac{\cU(t) \, \cU(t^2)}{2[1-\cU(t)][1-\cU(t^2)]}.
\end{align}
We have quickly obtained the same generating function reported in eq.~4 of Proposition 1 of \cite{AgranatTamirEtAl24b} and eq.~48 of \cite{AgranatTamirEtAl24a}.

The asymptotic growth of the coefficient $E_{n,1}$ of $\cE_1(t)$, describing the number of unlabeled time-consistent galled trees with $n$ leaves and 1 gall, follows \citep[Table 1]{AgranatTamirEtAl24b}:
\begin{equation}
E_{n,1} \sim \frac{1}{\gamma^3 \sqrt{\pi}} n^{1/2} \rho^{-n}.
\end{equation}

\begin{figure}[tb]
    \centering
    \includegraphics[width=9cm]{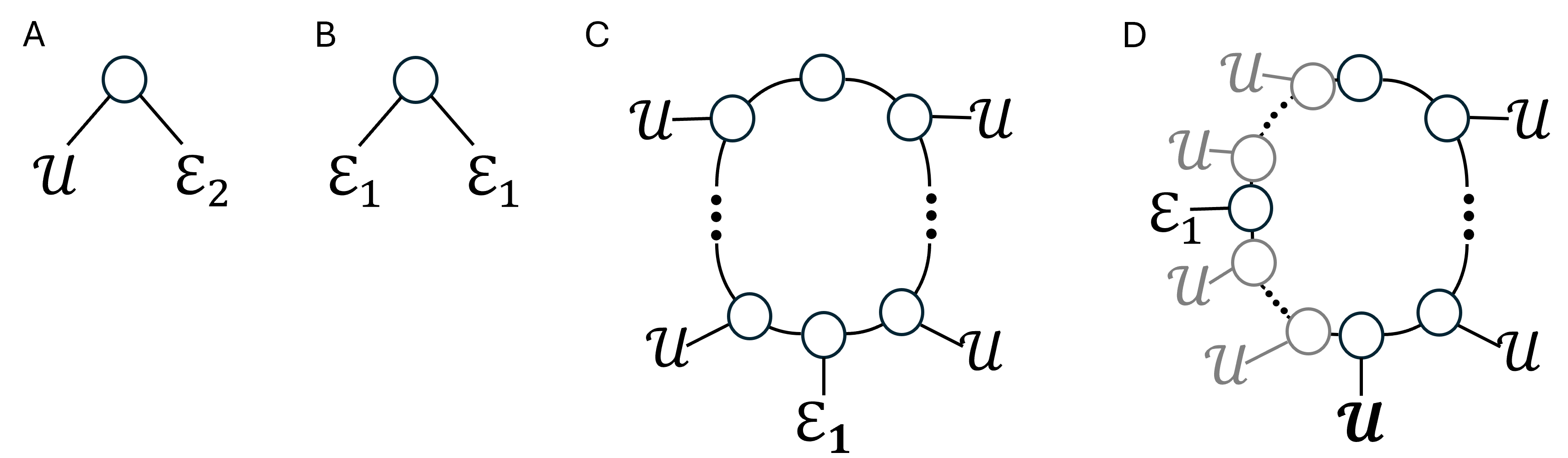}
    \caption{A non-plane unlabeled time-consistent galled tree with two galls, $\mathcal{E}_2$, has one of four structures. (A) A root with one subtree with no galls ($\mathcal{U}$) and one subtree with two galls ($\mathcal{E}_2$); there is no case of symmetry. (B) A root with two subtrees each with one gall ($\mathcal{E}_1$). This scenario has a case of symmetry, and because the tree is both non-plane and unlabeled, the two subtrees form a multiset of size 2. (C) A root gall and a subtree with one gall descended from the reticulation node (bold). On both sides of the reticulation node, because of the time-consistency condition, there are non-empty sequences of subtrees with no galls; because of symmetry, the two non-empty sequences form a multiset of size 2. (D) A root gall and a subtree with no galls descended from the reticulation node (bold). On one side of the reticulation node, there is a non-empty sequence of subtrees with no galls. On the other side, there is a subtree with one gall (to complete the tally of two galls), before and after which are are two possibly empty (because the subtree with the one gall is sufficient for time-consistency) sequences of subtrees with no galls (gray). There is no case of symmetry.}
    \label{fig:E2}
\end{figure}

\subsubsection{Two galls}

Next, we turn to $\cE_2(t)$, the generating function for time-consistent galled trees with exactly two galls. A galled tree with two galls falls into one of the following cases (Figure \ref{fig:E2}): 
\begin{enumerate}
 \item It has no root gall, and one of the two subtrees of the root contains both galls (Figure \ref{fig:E2}A); 
 \item It has no root gall, and the two subtrees of the root each contain one gall (Figure \ref{fig:E2}B);
 \item It has a root gall with two non-empty paths of unlabeled trees, and the subtree descended from the reticulation node contains the second gall (Figure \ref{fig:E2}C);
 \item It has a root gall, and the second gall is in a galled tree attached to one of the two paths (Figure \ref{fig:E2}D).
\end{enumerate}

Translating these cases into generating functions according to the symbolic method gives
\footnotesize
\begin{align}
    \mathcal{E}_2 = & \{ \circ \} \times 
    \bigg[
    \underbrace{ \mathcal{U} \times \mathcal{E}_2 }_{(1)}
    \textrm{   }\, \dot{\cup} \,\textrm{   }
    \underbrace{\textrm{MSET}_2 \Big( \mathcal{E}_1 \Big )}_{(2)}
     \textrm{   }\, \dot{\cup} \,\textrm{   }
      \underbrace{\mathcal{E}_1 \times \textrm{MSET}_2   \Big( \textrm{SEQ}^+(\mathcal{U})\Big)}_{(3)}
      \textrm{   }\, \dot{\cup} \,\textrm{   }
     \underbrace{\mathcal{U} \times \Big( \textrm{SEQ}(\mathcal{U}) \times \mathcal{E}_1 \times \textrm{SEQ}(\mathcal{U})\Big) \times \textrm{SEQ}^+(\mathcal{U})}_{(4)} 
    \bigg].
\end{align}
\normalsize
We have
\begin{align*}
\cE_2(t)&= \underbrace{ \cU(t) \, \cE_2(t) }_{(1)} + \underbrace{\rcp2\big[\cE_1(t)^2+\cE_1(t^2)\big]}_{(2)} 
 +\underbrace{\frac{\cE_1(t)}2\Bigg[\bigg(\frac{\cU(t)}{1-\cU(t)}\bigg)^2+\frac{\cU(t^2)}{1-\cU(t^2)}\Bigg]}_{(3)} \nonumber \\ 
 &\qquad + \underbrace{\cU(t) \cdot \frac{\cE_1(t)}{[1-\cU(t)]^2}\cdot\frac{\cU(t)}{1-\cU(t)}}_{(4)}. 
 \end{align*}
This expression for $\cE_2(t)$ simplifies, agreeing with eq.~14 in Proposition 4 of \cite{AgranatTamirEtAl24b}:
\begin{align}
 \cE_2(t)&=\rcp{2\big[1-\cU(t)\big]}\big[\cE_1(t)^2+\cE_1(t^2)\big] 
 +\frac{\cU(t)^2 \,\cE_1(t)}{2[1-\cU(t)]^3}+\frac{\cU(t^2) \, \cE_1(t)}{2[1-\cU(t)][1-\cU(t^2)]}
 +\frac{\cE_1(t) \, \cU(t)^2}{[1-\cU(t)]^4}.
\end{align}

The asymptotic growth of the coefficient $E_{n,2}$ of $\cE_2(t)$ is \citep[Table 1]{AgranatTamirEtAl24b} 
\begin{equation}
E_{n,2} \sim \frac{1}{3\gamma^7 \sqrt{\pi}} n^{5/2} \rho^{-n}.
\end{equation}

\subsubsection{Arbitrary numbers of galls}

We include the generating function for unlabeled time-consistent galled trees with no restriction on the number of galls. This generating function appears implicitly in Section 5.1 of \cite{AgranatTamirEtAl24a} and Section 6.3 of \cite{FuchsAndGittenberger24}. An arbitrary galled tree has three possibilities:
\begin{enumerate}
    \item It is a tree with one leaf.
    \item It has no root gall and two time-consistent galled subtrees.
    \item It has a root gall with two non-empty sequences of time-consistent galled trees, and a time-consistent galled subtree descended from the reticulation node.
\end{enumerate}
The symbolic method gives
\begin{equation}
   \mathcal{A} = \underbrace{\{ \square \}}_{(1)} \textrm{   }\, {\dot{\cup}} \,\textrm{   }
   \{ \circ \} \times \bigg[
   \underbrace{\textrm{MSET}_2(\mathcal{A})}_{(2)} \textrm{   }\, \dot{\cup} \,\textrm{   }
   \underbrace{\mathcal{A} \times \textrm{MSET}_2\left(\textrm{SEQ}^+(\mathcal{A})\right)}_{(3)}
   \bigg]. 
\end{equation}
The generating function is 
\begin{equation}
\mathcal{A}(t) 
= \underbrace{t}_{(1)}+ \underbrace{\frac{1}{2}\big[\mathcal{A}(t)^2+\mathcal{A}(t^2) \big]}_{(2)}
+\underbrace{\frac{\mathcal{A}(t)}{2}\Bigg[ \bigg( \frac{\mathcal{A}(t)}{1-\mathcal{A}(t)} \bigg)^2 + \frac{\mathcal{A}(t^2)}{1-\mathcal{A}(t^2)} \Bigg]}_{(3)}
\end{equation}

The asymptotic growth of the coefficient $A_n$ of $\mathcal{A}(t)$ is \citep[eq.~42]{AgranatTamirEtAl24a} 
\begin{equation}
A_n \sim (0.0779...)n^{-3/2} (0.2073...)^{-n}.
\end{equation}

\begin{figure}[tbp]
    \centering
    \includegraphics[width=6.3cm]{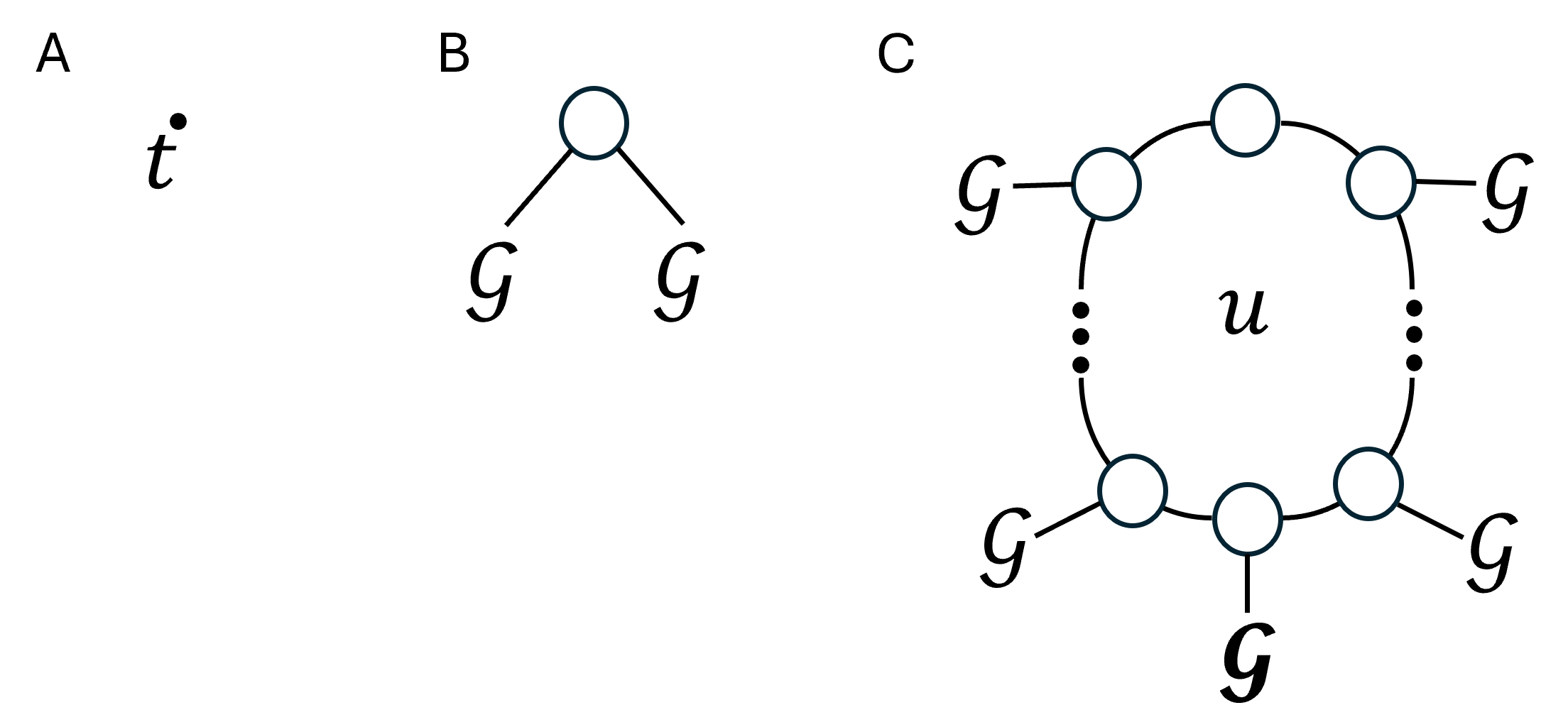}
    \caption{A non-plane unlabeled time-consistent galled tree with any number of galls, $\mathcal{G}$, has one of three structures. (A) A single leaf, $t$. (B) A root whose two subtrees are both time-consistent non-plane unlabeled galled trees with any number of galls. This scenario has a case of symmetry, and the two subtrees form a multiset of size 2. (C) A root gall, $u$, and a non-plane unlabeled time-consistent galled tree with any number of galls descended from the reticulation node (bold), following two non-empty (because of the time-consistency) sequences of non-plane unlabeled time-consistent  galled trees with any number of galls. This scenario has a case of symmetry, and hence, the two non-empty sequences form a multiset of size 2.}
    \label{fig:Gtu}
\end{figure}

\subsubsection{Arbitrary numbers of galls: the bivariate generating function} 

To find the generating function $\cE_g(t)$ for any fixed number of galls $g$, we use the bivariate generating function $\cG(t,u)=\sum_{n\ge 0}\sum_{g\ge0} E_{n,g} t^nu^g$. We derived this generating function in eq.~56 of \cite{AgranatTamirEtAl24a} and now show that it can be set up by the symbolic method similarly to the generating functions for the $g=1$ and $g=2$ cases.

A time-consistent galled tree structure with arbitrarily many galls, $\mathcal{G}$, has three cases. First, (1) it can be a single leaf (Figure \ref{fig:Gtu}A). Otherwise, (2) it has a binary root node with two galled trees attached (Figure \ref{fig:Gtu}B), or (3) it has a root gall (Figure \ref{fig:Gtu}C). In the third case, the construction needs a component $\mu$ to account for the root gall. We get
\begin{equation}
   \mathcal{G} = \underbrace{\{ \square \}}_{(1)} \textrm{   }\, {\dot{\cup}} \,\textrm{   }
   \{ \circ \} \times \bigg[
   \underbrace{\textrm{MSET}_2(\mathcal{G})}_{(2)} \textrm{   }\, \dot{\cup} \,\textrm{   }
   \underbrace{\mu \times \mathcal{G} \times \textrm{MSET}_2\left(\textrm{SEQ}^+(\mathcal{G})\right)}_{(3)}
   \bigg].
\end{equation}
The resulting generating function is
\begin{equation} \label{eq:G_fun}
\cG(t,u)=\underbrace{t}_{(1)} + \underbrace{ \rcp2\big[\cG(t,u)^2+\cG(t^2,u^2)\big]}_{(2)}
+\underbrace{\frac{u\cG(t,u)}{2}\Bigg[\bigg(\frac{\cG(t,u)}{1-\cG(t,u)}\bigg)^2+\frac{\cG(t^2,u^2)}{1-\cG(t^2,u^2)}\Bigg]}_{(3)}.
\end{equation} 

\subsubsection{Fixed number of galls}

We obtain the generating function $\cE_g(t)$ for a fixed number of galls $g$ by noting $\cU(t)=\cG(t,0)$, and for $g \geq 1$,
\[
\cE_g(t)=\rcp{g!}\bigg(\pdiff{}{u}{g}\cG\bigg)(t,0).
\]
Starting from $\cG(t,u)$, differentiating $g$ times with respect to $u$, setting $u=0$, and dividing by $g!$ yields an expression for $\cE_g(t)$ in terms of $\cU(t)$ and $\cE_i(t)$ with $i=1,2,\ldots,g-1$. The derivation relies on Leibniz's general product rule for higher-order derivatives and Faà di Bruno's formula for derivatives of a composition. We recall the latter formula, writing $\Du$ for $\partial/\partial u$: 
\begin{equation}
\Du^m(f\circ h)=\sum_{k_1+2k_2+\cdots +mk_m=m} \frac{m!}{k_1! \, k_2! \, \cdots \, k_m!} \bigg[ D_x^{k_1+k_2+\ldots
+k_m}\Big(f(x)\Big)\Big|_{x=h(u)} \bigg] \bigg[ \prod_{\ell=1}^m \bigg(\frac{\Du^{\ell} h}{\ell!}\bigg)^{k_\ell}\bigg].
\label{eq:fdb}
\end{equation}

To compute $\rcp{g!} \big(\Du^g\cG\big)(t,0)$, we begin from (\ref{eq:G_fun}). Via Leibniz's product rule, the first term yields
\[
\rcp{2g!}\Du^g\big(\cG(t,u)^2\big)=\rcp{2g!}
\sum_{\ell=0}^g \binom{g}{\ell} \Big(\Du^\ell\cG(t,u)\Big)\Big(\Du^{g-\ell}\cG(t,u)\Big).
\]
Setting $u=0$ gives 
\begin{equation} \label{eq:first_term}
\rcp{2g!}\Du^g\big(\cG(t,u)^2\big)\Big\vert_{u=0}=\rcp2
\sum_{\ell=0}^g \cE_\ell(t) \, \cE_{g-\ell}(t).
\end{equation} 

Next, we derive $\cG(t^2,u^2)$ with Faà di Bruno's formula (\ref{eq:fdb}), setting $f(u)=\cG(t^2,u)$ and $h(u)=u^2$, with $g$ in the role of $m$. As we eventually set $u=0$, the product in (\ref{eq:fdb}) can be nonzero only if the only factor that appears is $\Du^2h$. Hence, we must have $k_1=k_3=\dots=k_g=0$ and $k_2=g/2$, so that $g$ is even. We get 
\begin{equation} 
\label{second_term}
\rcp{2g!}\Du^g(\cG(t^2,u^2))\Big\vert_{u=0}=\frac{1}{2g!} \frac{g!}{\big(\frac g2\big)!} \cE_{\frac g2}(t^2) \, \Big(\frac g2\Big)!= 
\rcp2\cE_{\frac g2}(t^2).
\end{equation} 
for even $g$; for odd $g$, this derivative is zero.

The next term uses Leibniz's product rule and (\ref{eq:fdb}) with $f(u)=u^3/(1-u)^2$ and $h(u)=\cG(t,u)$:
\begin{align} 
\rcp{g!} \left.\Du^g\Bigg(\frac u2 \frac{\cG(t,u)^3}{\big(1-\cG(t,u)\big)^2}\Bigg)\right|_{u=0}&=
\rcp{2(g-1)!} \left.\Du^{g-1}\Bigg(\frac{\cG(t,u)^3}{\big(1-\cG(t,u)\big)^2}\Bigg)\right|_{u=0} \nonumber \\
&=\rcp2\sum_{k_1+2k_2+\cdots +(g-1)k_{g-1}=g-1} \rcp{k_1!\, k_2!\cdots k_{g-1}!} 
 \nonumber \\ 
&\quad\times \left.\bigg(\Du^{k_1+k_2+\ldots+k_{g-1}}\frac{u^3}{(1-u)^2}\bigg)\right|_{u=\cU(t)} \prod_{\ell=1}^{g-1} \cE_\ell(t)^{k_\ell} \label{eq:third_term}.
\end{align} 
We then insert 
\begin{equation} \label{eq:third_term_aux}
\rcp{k!}\Du^k \bigg[\frac{u^3}{(1-u)^2}\bigg] = 
\begin{cases}
\frac{u^2(3-u)}{(1-u)^3}= \frac{3u-1}{(1-u)^3}+1 & \text{ if } k=1,\\
\frac{3u+k-2}{(1-u)^{k+2}} & \text{ if } k \geq 2.
\end{cases}
\end{equation} 

What remains is the last term of (\ref{eq:G_fun}). First, note that 
\[
\left.\Du^k \big(u\cG(t,u) \big)\right|_{u=0}=\left.\Big(u\Du^k\cG(t,u)+k\Du^{k-1}\cG(t,u)\Big)\right|_{u=0}
=k! \, \cE_{k-1}(t). 
\]
Leibniz's rule then gives 
\begin{align} \label{fourth_term_1} 
\rcp{g!} \left.\Du^g\bigg(\frac{u\cG(t,u)}2 \frac{\cG(t^2,u^2)}{1-\cG(t^2,u^2)}\bigg)\right|_{u=0} &=
\bigg[\rcp{2}\sum_{k=1}^{g-1} \rcp{(g-k)!}\cE_{k-1}(t) \left.\Du^{g-k}\bigg(\frac{\cG(t^2,u^2)}{1-\cG(t^2,u^2)}\bigg)\right|_{u=0} \bigg] \nonumber \\
& \quad +\rcp2 \cE_{g-1}(t)\frac{\cU(t^2)}{1-\cU(t^2)} .
\end{align} 
For the derivative of the fraction, we again use Faà di Bruno's formula (\ref{eq:fdb}), now with $f(u)=u/(1-u)$
and $h(u)=\cG(t^2,u^2)$, and obtain, for $1 \leq k<g$,
\begin{align} 
\left.\Du^{g-k}\bigg(\frac{\cG(t^2,u^2)}{1-\cG(t^2,u^2)}\bigg)\right|_{u=0}&=
\sum_{r_1+2r_2+\ldots +(g-k)r_{g-k}=g-k} \frac{(g-k)!}{r_1! \, r_2! \cdots r_{g-k}!} 
\left.\bigg(\Du^{r_1+\dots+r_{g-k}}\frac{u}{1-u}\bigg)\right|_{u=\cU(t^2)} \nonumber \\
&\quad\times \prod_{m=1}^{g-k} \Bigg(\frac{\Du^m \cG(t^2,u^2)}{m!}\Bigg)^{r_m}\Bigg|_{u=0} 
 \label{fourth_term_2}.
\end{align}

Eq.~(\ref{second_term}) implies that in the product, terms with odd $m$ vanish unless $r_m=0$. So, we get a contribution only if all terms with odd indices among the $r_i$ are zero. In particular, $g-k$ is even as well.
 
Eq.~(\ref{fourth_term_2}) becomes
\begin{align} 
\left.\Du^{2b}\bigg(\frac{\cG(t^2,u^2)}{1-\cG(t^2,u^2)}\bigg)\right|_{u=0}&=
\sum_{2r_2+4r_4\cdots +2br_{2b}=2b} \frac{(2b)!}{r_1! \, r_2! \cdots r_{2b}!} 
\left.\bigg(\Du^{r_1+\dots+r_{2b}}\frac{u}{1-u}\bigg)\right|_{u=\cU(t^2)} \nonumber \\
&\quad\times \prod_{m=1}^{b} \Bigg(\frac{\Du^{2m} \cG(t^2,u^2)}{(2m)!}\Bigg)^{r_{2m}}\Bigg|_{u=0} 
 \label{fourth_term_2b}.
\end{align}

Finally, we use 
\begin{align}
\Du^k \bigg(\frac{u}{1-u}\bigg) &=\frac{k!}{(1-u)^{k+1}} \label{aux_1} \\
\Du^m(\cG(t^2,u^2))\Big|_{u=0} &=m! \, \cE_{\frac{m}{2}}(t^2) \textrm{ for even }m, \label{aux_2}
\end{align}

where (\ref{aux_2}) is equivalent to (\ref{second_term}). 
We collect (\ref{eq:first_term}), (\ref{second_term}), (\ref{eq:third_term}) (inserting (\ref{eq:third_term_aux})), and (\ref{fourth_term_1}) (inserting (\ref{fourth_term_2}) and (\ref{fourth_term_2b})) and get the expression for $\cE_g(t)$ after all: 
\small
\begin{align}
\cE_g(t) &= \bigg( \rcp2\sum_{\ell=0}^g \cE_\ell(t) \, \cE_{g-\ell}(t) \bigg) +\rcp2\cE_{\frac g2}(t^2) \label{eq:eg1}\\
& \quad +\rcp2\sum_{k_1+2k_2+\ldots +(g-1)k_{g-1}=g-1} \binom{k_1+k_2+ \ldots +k_{g-1}}{k_1,k_2,\ldots,k_{g-1}} 
\Bigg(\frac{3\cU(t)+(\sum_{i=1}^{g-1}k_i)-2}{\big(1-\cU(t)\big)^{(\sum_{i=1}^{g-1}k_i)+2}}+\delta_{1,\,\sum_{i=1}^{g-1}k_i}\Bigg) \nonumber \\
& \qquad \times \prod_{m=1}^{g-1} \cE_m(t)^{k_m} \label{eq:eg2}\\ 
& \quad +\rcp2\sum_{b=1}^{\big\lfloor\frac{ g-1}{2} \big\rfloor} \cE_{(g-2b)-1}(t) \sum_{2r_2+4r_4\cdots +(2b)r_{2b}=2b} \binom{r_1+r_2+\ldots +r_{2b}}{r_1,r_2,\ldots,r_{2b}} 
\rcp{\big(1-\cU(t^2)\big)^{(\sum_{i=1}^{b}r_{2i})+1}} \nonumber \\
& \quad \times \prod_{m=1}^{ b}
\cE_m(t^2)^{r_{2m}} \label{eq:eg3a} \\ 
& \qquad +\rcp2 \cE_{g-1}(t)\frac{\cU(t^2)}{1-\cU(t^2)} \label{eq:eg3b}.
\end{align}
\normalsize

\bp \label{prop:equal}
The equation for the generating function $\cE_g(t)$ described by (\ref{eq:eg1})--(\ref{eq:eg3b}) is equal to the equation for the generating function in eq.~17 of \cite{AgranatTamirEtAl24b}.
\ep

\bpf
Eq.~17 of \cite{AgranatTamirEtAl24b} is written as half of the sum of three equations, eqs.~20, 21, and 22 of \cite{AgranatTamirEtAl24b}. We will show three properties:
\begin{enumerate}
    \item (\ref{eq:eg1}) is equal to half of eq.~20 in \cite{AgranatTamirEtAl24b}.
    \item (\ref{eq:eg2}) is equal to half of eq.~21 in \cite{AgranatTamirEtAl24b}.
    \item The sum of (\ref{eq:eg3a}) and (\ref{eq:eg3b}) is equal to half of eq.~22 in \cite{AgranatTamirEtAl24b}.
\end{enumerate}
\noindent The first property (1) is trivial, remembering that $\mathcal{E}_0(t) = \mathcal{U}(t)$. For (2), we must show three further points: \\
(2i) The notation $\ell$ of \cite{AgranatTamirEtAl24b} is equal to $\sum_{i=1}^{g-1}k_i$ here. \\
(2ii) The limits of summation $\sum_{\ell=1}^{g-1}\sum_{d\in C(g-1,\ell)}$ in eq.~21 of \cite{AgranatTamirEtAl24b} are the same as $\sum_{k_1+2k_2+\cdots +(g-1)k_{g-1}=g-1} \binom{k_1+k_2+\ldots +k_{g-1}}{k_1,k_2,\ldots,k_{g-1}} $ from (\ref{eq:eg2}) here. \\
(2iii) In the notation of \cite{AgranatTamirEtAl24b}, $\prod_{j=1}^\ell\mathcal{E}_{d_j}(t)$ is equal to $\prod_{m=1}^{g-1} \mathcal{E}_m(t)^{k_m}$ here.

\medskip 
The sum $\sum_{\ell=1}^{g-1}\sum_{d\in C(g-1,\ell)}$ from \cite{AgranatTamirEtAl24b}, where $C(g-1,\ell)$ represents the set of compositions of $g-1$ into $\ell$ positive parts, traverses possible numbers of subtrees $\ell$ from the root gall that have at least one gall, and all ways to distribute $g-1$ galls among them. The sum $\sum_{k_1+2k_2+\cdots +(g-1)k_{g-1}=g-1}$ from (\ref{eq:eg2}) traverses the numbers of subtrees $(k_1,k_2, \ldots, k_{g-1})$ with all possible positive numbers of galls $(1,2, \ldots, g-1)$ and total $g-1$. The total number of trees with a positive number of galls is then $\sum_{i=1}^{g-1}k_i$, the same quantity as $\ell$ from \cite{AgranatTamirEtAl24b}. Hence, (2i) is proven; note that the Iverson bracket $[\![\ell=1]\!]$ of \cite{AgranatTamirEtAl24b} matches the Kronecker delta $\delta_{1,\sum_{i=1}^{g-1}k_i}$ here.

For (2ii), with the number of galls $\sum_{i=1}^{g-1}k_i$ and the specific $k_i$ values specified, the number of ways to distribute the total number of galls into the $\sum_{i=1}^{g-1}k_i$ subtrees of the root gall that contain at least one gall is $\binom{k_1+k_2+\ldots +k_{g-1}}{k_1,k_2,\ldots,k_{g-1}}$. In total, $\sum_{k_1+2k_2+\cdots +(g-1)k_{g-1}=g-1} \binom{k_1+k_2+\ldots +k_{g-1}}{k_1,k_2,\ldots,k_{g-1}}$ traverses the same arrangements of galls into subtrees as $\sum_{\ell=1}^{g-1}\sum_{d\in C(g-1,\ell)}$ from \cite{AgranatTamirEtAl24b}.

Next, for (2iii), in the notation of \cite{AgranatTamirEtAl24b}, the number of subtrees $\ell$ with a positive number of galls and the numbers of galls $\{ d_j \}_{j=1}^{\ell}$ in each of these subtrees are determined, $\prod_{j=1}^\ell\mathcal{E}_{d_j}(t)$ traverses these subtrees and takes their associated product. The notation $\prod_{m=1}^{g-1} \mathcal{E}_m(t)^{k_m}$ here
takes the same product for $(k_1, k_2, \ldots, k_{g-1})$ fixed, proceeding in a different order by traversing each possible number of galls in the subtrees ($m=1,2,\ldots,g-1$), counting how many subtrees have that number of galls ($k_m$).

Finally, for (3), we first show that eq.~22 of \cite{AgranatTamirEtAl24b} for $\ell \neq 0$ is equal to (\ref{eq:eg3a}) here. We begin by showing that $\sum_{\ell=1}^{\lfloor \frac{g-1}{2} \rfloor}\sum_{b=\ell}^{\lfloor \frac{g-1}{2} \rfloor}\sum_{\mathbf{d}\in C(b,\ell)}$ in the notation of \cite{AgranatTamirEtAl24b} is equal to $\sum_{b=1}^{\lfloor\frac{ g-1}{2} \rfloor}\sum_{2r_2+4r_4\ldots +(2b)r_{2b}=2b}$ in (\ref{eq:eg3a}). The notation of \cite{AgranatTamirEtAl24b} traverses possible numbers of subtrees $\ell$ from the root gall on one side of the reticulation node, with a positive number of galls, then determines the number of galls $b\geq \ell$ to distribute among those subtrees, and finally distributes the galls with $\sum_{\mathbf{d}\in C(b,\ell)}$. (\ref{eq:eg3a}) does the same computation in a different order, first traversing the number of galls $b$ on one side of the reticulation node and then how they are distributed into at most $b$ subtrees $\sum_{2r_2+4r_4\cdots +(2b)r_{2b}=2b}$. Hence, $\ell$ in eq.~22 of \cite{AgranatTamirEtAl24a} is equal to ${r_2+r_4\ldots +r_{2b}}$ in (\ref{eq:eg3a}) and $\prod_{j=1}^{\ell}\mathcal{E}_{d_j}(t^2)$ in \cite{AgranatTamirEtAl24b} is equal to $\prod_{m=1}^{b}\mathcal{E}_m(t^2)^{r_{2m}}$ here; this equivalence is similar to that seen in the proof of (2iii), as the $d_j$ are necessarily positive and their sum is $b$, as is the sum $r_2+2r_4\cdots +br_{2b}$. Finally, $\mathcal{E}_{g-2b-1}(t)$ appears in both representations, for each $b$ with $b=1,2,\ldots,\lfloor \frac{g-1}{2} \rfloor$.

\begin{table}[tbh]
    \centering
    \small
    \begin{tabular}{|c|r|rrrrr|}
    \hline
          &              & \multicolumn{5}{c|}{Number of trees} \\ 
Number of & Total number & \multicolumn{5}{c|}{with a fixed number of galls ($E_{n,g}$)} \\ 
leaves ($n$) & \multicolumn{1}{c|}{of trees ($A_n$)} & $g=0$ & $g=1$ & $g=2$ & $g=3$ & $g=4$ \\ \hline
 1 &      1 &  1 &    - &    - &     - &   - \\
 2 &      1 &  1 &    - &    - &     - &   - \\
 3 &      2 &  1 &    1 &    - &     - &   - \\
 4 &      6 &  2 &    4 &    - &     - &   - \\
 5 &     72 &  3 &   15 &    2 &     - &   - \\
 6 &    272 &  6 &   48 &   18 &     - &   - \\
 7 &   1064 & 11 &  148 &  107 &     6 &   - \\
 8 &   4271 & 23 &  435 &  528 &    78 &   - \\
 9 & 17,497 & 46 & 1250 & 2295 &   661 &  19 \\
10 & 72,483 & 98 & 3512 & 9185 &  4356 & 346 \\
\hline
\end{tabular}
\caption{Numbers of unlabeled time-consistent galled trees with specified numbers of leaves and galls. Entries $E_{n,g}$ are computed recursively from (\ref{eq:unlabeled-rec}) and are copied from Table 3 of \cite{AgranatTamirEtAl24a}.}
\label{table:unlabeled-gallnum}
\end{table}

It is left to show that the $\ell=0$ case in eq.~22 of \cite{AgranatTamirEtAl24a} is equal to (\ref{eq:eg3b}). When $\ell=0$, $\sum_{\mathbf{d}\in C(b,\ell)}$ is not zero only if $b=\ell=0$ and so eq.~22 of \cite{AgranatTamirEtAl24a} is equal to 
\begin{align*}
    \bigg(\frac{1}{1-\cU(t^2)}-1 \bigg)\cE_{g-1}&= \frac{\cE_{g-1}(t) \, \cU(t^2)}{1-\mathcal{U}(t^2)},
\end{align*}
as is needed. \epf

The asymptotic growth of the coefficient $E_{n,g}$ of $\cE_g(t)$ is \citep[eq.~42]{AgranatTamirEtAl24b} 
\begin{equation}
E_{n,g} \sim \frac{2^{2g-1}}{(2g)! \, \gamma^{4g-1} \sqrt{\pi} } n^{2g-\frac{3}{2} } \rho^{-n}. 
\end{equation}

\subsection{Numerical computation}

Proposition 3 of \cite{AgranatTamirEtAl24b} provides a recursive computation by which numerical values of $E_{n,g}$, the number of unlabeled time-consistent galled trees with $n$ leaves and $g$ galls. We begin from $E_{1,0}=1$ and $E_{1,g}=0$ for $g \geq 1$. Let $C(n,k)$ denote the compositions of $n$ into $k$ (positive) parts. Let $C_p(n,k)$ denote the palindromic compositions of $n$ into $k$ parts, where a palindromic composition is a composition that is invariant when written in reverse order. We have 
\begin{align}
E_{n,g} &= \frac{1}{2} \bigg[ \bigg( \sum_{\mathbf{c} \in C(n,2)} \sum_{\mathbf{d} \in C(g+2,2)} \prod_{i=1}^2 E_{c_i,d_i-1} \bigg) + E_{\frac{n}{2},\frac{g}{2}} \nonumber \\
& \quad + \bigg( \sum_{k=3}^n (k-2) \sum_{\mathbf{c} \in C(n,k)} \sum_{\mathbf{d} \in C(g-1+k,k)} \prod_{i=1}^k E_{c_i,d_i-1} \bigg) \nonumber \\
& \quad + \bigg( \sum_{a=1}^{\lfloor \frac{n-1}{2} \rfloor} \sum_{\mathbf{c} \in C_p(n,2a+1)} \sum_{\mathbf{d} \in C_p\big(g-1+(2a+1),2a+1 \big)} \prod_{i=1}^{a+1} E_{c_i,d_i-1}
\bigg) \bigg].
\label{eq:unlabeled-rec}
\end{align}
Iterating the recursion, Table \ref{table:unlabeled-gallnum} follows Table 3 of \cite{AgranatTamirEtAl24a} and gives the numbers of unlabeled time-consistent galled trees for small $n$ and $g$.

\section{Symbolic method for labeled time-consistent galled trees}
\label{sec:labeled}

\subsection{Overview}

\cite{AgranatTamirEtAl24a}, \cite{AgranatTamirEtAl24b}, and Section \ref{sec:unlabeled} have focused on enumeration of unlabeled time-consistent galled trees; in this section, we perform analogous enumerations, but now with the leaves labeled. As each unlabeled shape has many possible labelings, the labeled structures are more numerous. However, in conducting the enumeration, the number of cases is smaller, so that the enumeration can be performed in fewer steps. We obtain the generating functions and also the asymptotic enumerations.

\subsection{Definitions}

We consider the leaf-labeled time-consistent galled trees. This class of networks is obtained by considering all possible labelings of the unlabeled time-consistent galled trees of Section \ref{sec:unlabeled-definitions}. Similarly to the unlabeled case, shapes are \emph{non-plane}, and the \emph{time-consistent} leaf-labeled galled trees are identical to the \emph{normal} leaf-labeled galled trees. The main difference of the enumeration of labeled time-consistent galled trees from the enumeration of unlabeled time-consistent galled trees is that the use of labels eliminates the symmetry case.

\subsection{Symbolic method}

We apply the symbolic method for labeled structures. The approach proceeds similarly to the case of unlabeled structures, with a difference that in translating constructions into algebraic expressions, we consider exponential generating functions~\citep[][pp.~97-98]{FlajoletAndSedgewick09}. 

\subsection{Generating functions}

\subsubsection{No galls}

We begin with labeled galled trees with no galls. These are \emph{labeled topologies}, or \emph{labeled non-plane trees}. Denote by $u_n$ the number of labeled trees with $n$ leaves. It is well known that $u_n=(2n-3)!! = (2n-2)!/[ 2^{n-1} (n-1)! ]$ for $n\geq2$, with $u_1=1$ \citep{EdwardsAndCSforza64, Felsenstein78}. 

A labeled tree with no galls $\mathfrak{U}$ is either a single leaf ($\square$), or it is a root $(\circ)$ to which a set of two distinct labeled trees are attached, $\textrm{SET}_2(\mathfrak{U})$. It can be written in the symbolic method as 
\begin{equation}
   \mathfrak{U} = \{ \square \} \textrm{   }\, {\dot{\cup}}  \,\textrm{   }
   \{ \circ \} \times
   \textrm{SET}_2(\mathfrak{U}),
\end{equation}
The exponential generating function satisfies
\begin{equation}
    \mathfrak{U}(t) = \sum_{n=0}^\infty \frac{u_n}{n!}t^n = t + \frac{1}{2}\mathfrak{U}(t)^2.
\end{equation}
The factor of $\frac{1}{2}$ arises because the trees are non-plane. The exponential generating function is 
\begin{equation}
\label{eq:U-labeled}
\mathfrak{U}(t) = 1-\sqrt{1-2t},
\end{equation}
so that $\frac{1}{2}$ is the radius of convergence. Using asymptotic theory of generating functions~\citep[][p.~392]{FlajoletAndSedgewick09} to obtain the asymptotic approximation to $u_n$, as $n \rightarrow \infty$, we apply Figure VI.5 of \cite{FlajoletAndSedgewick09}, producing
\begin{equation}
u_n \sim \frac{n^{-3/2}}{2\sqrt{\pi}}\Big(\frac{1}{2}\Big)^{-n}n!\sim \frac{\sqrt{2}}{2}n^{n-1}\Big(\frac{2}{e}\Big)^n,
\end{equation}
where the last approximation follows by Stirling's approximation to $n!$.

\subsubsection{One gall}

As with unlabeled trees, the next case is galled trees with exactly one gall. Following Figure \ref{fig:E1}, a galled tree with one gall is either (1) a root with one subtree with no galls and one subtree with exactly one gall, or (2) a root gall with a subtree with no galls descended from the reticulation node and a non-empty sequence of subtrees with no galls on each side of the gall. By the symbolic method, denoting by $\mathfrak{E}_1(t)$ the exponential generating function of labeled time-consistent galled trees with exactly one gall, we have:
\begin{equation}
    \mathfrak{E}_1 = \{ \circ \} \times 
    \Big[
    \underbrace{ \mathfrak{U} \star \mathfrak{E}_1 }_{(1)} 
    \textrm{  } \dot{\cup} \textrm{  }
    \underbrace{\mathfrak{U} \star \textrm{SET}_2\Big(\textrm{SEQ}^+(\mathfrak{U})\Big)}_{(2)}
    \Big]
\end{equation}
where $\star$ is the labeled product.

Denoting by $e_{n,1}$ the number of labeled time-consistent galled trees with $n$ leaves and exactly one gall, we convert to a generating function using the library~\citep[Figure II.18, p.~148]{FlajoletAndSedgewick09}, 
\begin{align}
    \mathfrak{E}_1(t) &= \sum_{n=0}^\infty \frac{e_{n,1}}{n!}t^n = \underbrace{ \mathfrak{U}(t) \, \mathfrak{E}_1(t) }_{(1)}
    +\underbrace{\mathfrak{U}(t) \, \frac{1}{2}\bigg(\frac{\mathfrak{U}(t)}{1-\mathfrak{U}(t)}\bigg)^2}_{(2)} \nonumber \\ 
     &= \frac{\mathfrak{U}(t)^3}{2\big(1-\mathfrak{U}(t)\big)^3}.
\end{align}
Again, the $\frac{1}{2}$ arises to avoid double-counting structures that arise with the two sides of the root gall exchanged.

Applying (\ref{eq:U-labeled}), we have:
\begin{equation}
    \mathfrak{E}_1(t) = \frac{(1-\sqrt{1-2t})^3}{2(1-2t)^{3/2}},
\end{equation}
so that as $t \rightarrow \frac{1}{2}$, 
\begin{equation}
\label{eq:E1asy}
    \mathfrak{E}_1(t) \sim \frac{1}{2}(1-2t)^{-3/2}.
\end{equation}
Next, for the asymptotic growth of $e_{n,1}$ as $n \rightarrow \infty$, applying Figure VI.5 of \cite{FlajoletAndSedgewick09}
\begin{equation}
    e_{n,1} \sim \frac{1}{2}\cdot\frac{n^{1/2}}{\Gamma(\frac{3}{2})}\Big(\frac{1}{2}\Big)^{-n}n! = \frac{n^{1/2}}{\sqrt{\pi}}\Big(\frac{1}{2}\Big)^{-n}n!
    \sim \sqrt{2}n^{n+1}\Big(\frac{2}{e}\Big)^n. \label{eq:en1asy}
\end{equation}
where the last step uses the Stirling approximation.

\subsubsection{Two galls}

Continuing to mimic constructions in the unlabeled case, a galled tree with two galls has four possibilities, as represented by Figure \ref{fig:E2}: (1) a root with a subtree with no galls and a subtree with two galls; (2) a root with two subtrees each with one gall; (3) a root gall with a subtree with one gall from the reticulation node and two non-empty sequences of trees with no galls on either side; and (4) a root gall with a tree with no galls from the reticulation node and a non-empty set of trees with no galls on one side, and, on the other side, a tree with one gall and two (perhaps empty) sequences of trees with no galls on each side of that tree. By the symbolic method:
\begin{align}
    \mathfrak{E}_2 &= \{ \circ \}\times\\ \nonumber
    & \qquad \Big[
    \underbrace{\mathfrak{U} \star \mathfrak{E}_2 }_{(1)}
    \textrm{  } \dot{\cup} \textrm{  }
    \underbrace{\textrm{SET}_2(\mathfrak{E}_1)}_{(2)}
    \textrm{  } \dot{\cup} \textrm{  }
    \underbrace{\mathfrak{E}_1 \star \textrm{SET}_2\Big(\textrm{SEQ}^+(\mathfrak{U})\Big)}_{(3)}
    \textrm{  } \dot{\cup} \textrm{  }
    \underbrace{\mathfrak{U} \star \Big( \textrm{SEQ}(\mathfrak{U}) \star\mathfrak{E}_1\star \textrm{SEQ}(\mathfrak{U}) \Big)\star
    \textrm{SEQ}^+(\mathfrak{U})}_{(4)}
    \Big].
\end{align}

Denoting by $e_{n,2}$ the number of labeled galled trees with $n$ leaves and exactly two galls, 
\begin{align}
    \mathfrak{E}_2(t) &= \sum_{n=0}^\infty \frac{e_{n,2}}{n!}t^n \nonumber \\
    &= \underbrace{\mathfrak{U}(t) \, \mathfrak{E}_2(t) }_{(1)} 
    + \underbrace{\frac{1}{2}\mathfrak{E}_1(t)^2}_{(2)} 
    + \underbrace{\mathfrak{E}_1(t)\frac{1}{2}\bigg(\frac{\mathfrak{U}(t)}{1-\mathfrak{U}(t)}\bigg)^2}_{(3)}
    + \underbrace{\mathfrak{U}(t)\cdot\frac{\mathfrak{E}_1(t)}{\big(1-\mathfrak{U}(t)\big)^2}\cdot\frac{\mathfrak{U}(t)}{1-\mathfrak{U}(t)}}_{(4)} \nonumber \\
    &= \frac{\mathfrak{E}_1(t)^2}{2\big(1-\mathfrak{U}(t)\big)}
    + \frac{\mathfrak{E}_1(t)\, \mathfrak{U}(t)^2}{2\big(1-\mathfrak{U}(t)\big)^3} 
    + \frac{\mathfrak{E}_1(t)\, \mathfrak{U}(t)^2}{\big(1-\mathfrak{U}(t)\big)^4}
\end{align}
Using approximations for $\mathfrak{U}(t)$ (\ref{eq:U-labeled}) and $\mathfrak{E}_1(t)$ (\ref{eq:E1asy}) as $t \rightarrow \frac{1}{2}$, we have 
\begin{align}
    \mathfrak{E}_2(t) &\sim \frac{\frac{1}{4}(1-2t)^{-3}}{2(1-2t)^{1/2}}
    + \frac{\frac{1}{2}(1-2t)^{-3/2}}{{2}(1-2t)^{3/2}} + \frac{\frac{1}{2}(1-2t)^{-3/2}}{(1-2t)^2} \nonumber \\
    &\sim \frac{5}{8}(1-2t)^{-7/2}. \label{eq:E2asy}
\end{align}

Therefore, as $n \rightarrow \infty$, by Corollary VI.1 of \cite{FlajoletAndSedgewick09},
\begin{equation}
    e_{n,2} \sim \frac{5}{8}\cdot\frac{n^{5/2}}{\Gamma(\frac{7}{2})}\Big(\frac{1}{2}\Big)^{-n}n! = \frac{n^{5/2}}{3\sqrt{\pi}}\Big(\frac{1}{2}\Big)^{-n}n! \sim \frac{\sqrt{2}}{3}n^{n+3}\Big(\frac{2}{e}\Big)^n. \label{eq:en2asy}
\end{equation}

\subsubsection{Arbitrary numbers of galls}

To find the exponential generating function for labeled time-consistent galled trees with any fixed number of galls, we first examine the exponential generating function for labeled time-consistent galled trees with no restriction on the number of galls. The generating function was obtained implicitly in Section 5.2 of~\cite{FuchsAndGittenberger24}. An arbitrary time-consistent galled tree has one of the following forms:
\begin{enumerate}
    \item It is a tree with one leaf.
    \item It has no root gall and two time-consistent galled subtrees.
    \item It has a root gall with two non-empty sequences of time-consistent galled trees, and a time-consistent galled subtree descended from the reticulation node.
\end{enumerate}
By the symbolic method,
\begin{equation}
   \mathfrak{A} = \underbrace{\{ \square \}}_{(1)} \textrm{   }\, {\dot{\cup}} \,\textrm{   }
   \{ \circ \} \times \bigg[
   \underbrace{\textrm{SET}_2(\mathfrak{A})}_{(2)} \textrm{   }\, \dot{\cup} \,\textrm{   }
   \underbrace{\mathfrak{A} \star \textrm{SET}_2\left(\textrm{SEQ}^+(\mathfrak{A})\right) }_{(3)}
   \bigg]. 
\end{equation}
The exponential generating function has the form:
\begin{equation}
    \mathfrak{A}(t) = \sum_{n=0}^\infty \frac{a_n}{n!}t^n = \underbrace{t}_{(1)} + \underbrace{\frac{1}{2}\mathfrak{A}(t)^2}_{(2)} 
    + \underbrace{\frac{1}{2}\mathfrak{A}(t)\bigg(\frac{\mathfrak{A}(t)}{1-\mathfrak{A}(t)}\bigg)^2}_{(3)}.
\end{equation}

Numerical computations of the number of labeled time-consistent galled trees with no restriction on the number of galls appeared in Theorem 8 of \cite{CardonaAndZhang20}. Asymptotic analysis of their growth was considered in Section 5.2 of \cite{FuchsAndGittenberger24}.

\subsubsection{Arbitrary numbers of galls: the bivariate generating function}

As in the case of unlabeled time-consistent galled trees, we derive the bivariate generating function for labeled time-consistent galled trees with a fixed number of galls by adding the root gall. We have
\begin{equation}
\mathfrak{G} = \underbrace{\{\square\}}_{(1)} \textrm{   } \dot{\cup} \textrm{   } 
    \{\circ\} \times
    \bigg[\underbrace{\textrm{SET}_2(\mathfrak{G})}_{(2)} \textrm{   } \dot{\cup} \textrm{   } \underbrace{\mu\times\mathfrak{G}\star \textrm{SET}_2\big(\textrm{SEQ}^+(\mathfrak{G})\big)}_{(3)}\bigg]
\end{equation}

Denoting by $e_{n,g}$ the number of labeled time-consistent galled trees with $n$ leaves and $g$ galls, we have
\begin{equation}
    \mathfrak{G}(t,u) = \sum_{n=0}^\infty\sum_{g=0}^\infty \frac{e_{n,g}}{n!}t^nu^g = t + \frac{1}{2}\mathfrak{G}(t,u)^2
    +\frac{u\mathfrak{G}(t,u)}{2}\bigg(\frac{\mathfrak{G
    }(t,u)}{1-\mathfrak{G}(t,u)}\bigg)^2.
\end{equation}

\subsubsection{Fixed number of galls}

To find the exponential generating function for labeled time-consistent galled trees with exactly $g$ galls, we again differentiate the bivariate function $g$ times with respect to $u$, set $u=0$, and divide by $g!$. The exponential generating function in the labeled case has three of the five terms summed in the generating function for the unlabeled case (\ref{eq:G_fun}); two terms that result from symmetries that exist in the unlabeled case but not the labeled case are omitted. We can simply copy the relevant terms from the analysis of the unlabeled case to obtain 
\small
\begin{align}
    \mathfrak{E}_g(t) &= \sum_{n=0}^\infty \frac{e_{n,g}}{n!}t^n \nonumber \\ 
    &= \label{eq:1termL}
    \frac{1}{2}\sum_{\ell=0}^g \mathfrak{E}_\ell(t) \, \mathfrak{E}_{g-\ell}(t) \\ 
& \quad \label{eq:2termL} +\rcp2\sum_{k_1+2k_2+\ldots +(g-1)k_{g-1}=g-1} \binom{k_1+k_2+\ldots +k_{g-1}}{k_1,k_2,\ldots,k_{g-1}} 
\Bigg(\frac{3\mathfrak{U}(t)+(\sum_{i=1}^{g-1}k_i)-2}{(1-\mathfrak{U}(t))^{(\sum_{i=1}^{g-1}k_i)+2}}+\delta_{1,\,\sum_{i=1}^{g-1}k_i}\Bigg) \nonumber \\
& \qquad \times \prod_{m=1}^{g-1} \mathfrak{E}_m(t)^{k_m} \\ 
\nonumber
&= \frac{1}{2\big(1-\mathfrak{U}(t)\big)}\bigg[\sum_{\ell=1}^{g-1} \mathfrak{E}_\ell(t) \, \mathfrak{E}_{g-\ell}(t) \\ \nonumber
& \quad + \sum_{k_1+2k_2+\cdots +(g-1)k_{g-1}=g-1} \binom{k_1+k_2+\ldots +k_{g-1}}{k_1,k_2,\ldots,k_{g-1}} \bigg(\frac{3\mathfrak{U}(t)+(\sum_{i=1}^{g-1}k_i)-2}{(1-\mathfrak{U}(t))^{(\sum_{i=1}^{g-1}k_i)+2}}+\delta_{1,\,\sum_{i=1}^{g-1}k_i}\Bigg)\prod_{\ell=1}^{g-1} \mathfrak{E}_\ell(t)^{k_\ell} \bigg].\\ \label{eq:EgtL}
\end{align}
\normalsize
Equation (\ref{eq:1termL}) results from copying (\ref{eq:first_term}), and (\ref{eq:2termL}) results from copying (\ref{eq:third_term}) and (\ref{eq:third_term_aux}).

To find the asymptotic approximation for $e_{n,g}$, we require a small proposition.
\bp \label{prop:equi}
The equation for the generating function $\mathfrak{E}_g(t)$ (\ref{eq:EgtL}) is equivalent to the sum of the terms in the generating function for unlabeled time-consistent galled trees that determine the asymptotic growth of the coefficients, namely eq.~21 of \cite{AgranatTamirEtAl24b} and the second term of eq.~20 of \cite{AgranatTamirEtAl24b}.
\ep
\bpf
Eq.~(\ref{eq:1termL}) is equivalent to the first term of (\ref{eq:eg1}), and (\ref{eq:2termL}) is equivalent to (\ref{eq:eg2}). We have proven in Proposition \ref{prop:equal} the equality of the first term of (\ref{eq:eg1}) to the second term of eq.~20 of \cite{AgranatTamirEtAl24b} and of (\ref{eq:eg2}) to eq.~21 of \cite{AgranatTamirEtAl24b}. In Section 6.2 of \cite{AgranatTamirEtAl24b}, it is these terms that determine the asymptotic growth of the coefficients of the generating function. 
\epf

Based on Proposition \ref{prop:equi}, we follow the analysis in Propositions 7 and 9 of \cite{AgranatTamirEtAl24b}. We replace the constants
$\gamma \sim 1.1300$ and $\rho \sim 0.4027$ that arise for the unlabeled $\mathcal{E}_1(t)$ with $\gamma=1$ and $\rho=\frac{1}{2}$ in the labeled case. As $t \rightarrow \frac{1}{2}$, we have
\begin{equation}
    \mathfrak{E}_g(t) \sim \frac{C_{2g-1}}{2^{2g-1}(1-2t)^{2g-1/2}} = \frac{(4g-3)!!}{(2g)! \, (1-2t)^{2g-1/2}},
\end{equation}
where $C_{2g-1}$ represents the Catalan number $\frac{1}{2g} {4g-2 \choose 2g-1}$. From Theorem 10 of \cite{AgranatTamirEtAl24b}, 
\begin{equation}
\label{eq:eng}
    e_{n,g} \sim \frac{2^{2g-1}}{(2g)!\sqrt{\pi}}n^{2g-3/2}\Big(\frac{1}{2}\Big)^{-n}n! \sim \frac{2^{2g-1}\sqrt{2}}{(2g)!}\Big(\frac{2}{e}\Big)^n n^{n+2g-1}.
\end{equation}

This expression for the asymptotic growth of the number of labeled time-consistent galled trees is quite similar to the asymptotic approximation for the growth of the number of general phylogenetic networks $PN$, reticulation-visible networks $RV$, galled networks $GN$, tree-child networks $TC$, normal networks $N$, and galled tree-child networks $GTC$, considering networks with $n$ labeled leaves and $k$ reticulations, or (\ref{eq:asymptotic}), $[{2^{k-1}\sqrt{2}}/{k!}] n^{n+2k-1} (2/e)^n n^{n+2k-1}$.

In galled trees, the number of galls $g$ corresponds to the number of reticulations. The expression (\ref{eq:eng}) accords with (\ref{eq:asymptotic}) for $g=0$, in which case both the time-consistent galled trees and the various network classes reduce to labeled binary trees. The expression (\ref{eq:eng}) also accords with (\ref{eq:asymptotic}) for $g=1$, for which the network classes describe labeled trees with one gall.

For two or more reticulations, however, the expressions differ in the subexponential growth. In particular, inserting $k$ in place of $g$ in (\ref{eq:eng}), time-consistent galled trees number fewer than the networks in the other classes by a factor of $2^k k! / (2k)! = 1/(2k-1)!!$. The time-consistent galled trees are asymptotically fewer in number than phylogenetic networks in each of several classes of more permissive structures.

\subsection{Numerical computation}

We have calculated the large-$n$ approximation to the number of labeled time-consistent galled trees with a fixed number of galls. We now calculate the \emph{exact} number of such trees for small $n$ and $g$. 

Theorem 8 of \cite{CardonaAndZhang20} found an exact formula for the number of labeled time-consistent galled trees with $n$ leaves, summing over all possible numbers of galls. We develop a recursive formula for the number of such trees with a \emph{given} number of galls $g$.

First, for a single leaf, $e_{1,0}=1$, and for all $g>0$, $e_{1,g}=0$. Otherwise, if a time-consistent galled tree has at least two leaves, then either it has two subtrees of the root or it has a root gall. In the former case, for each possible number of leaves $m$ that can be assigned to the ``left'' subtree of the root, the number of ways of choosing the $m$ labels for those leaves is ${n \choose m}$. We sum over all possible values $\ell$ for the number of galls in the ``left'' subtree, obtaining for the contribution of trees with no root gall
\begin{equation}
\label{eq:engrec1}
    \frac{1}{2}\bigg{[}\sum_{m=1}^{n-1} \binom{n}{m} \sum_{\ell=0}^g e_{m,\ell}e_{n-m,g-\ell}\bigg{]}.
\end{equation}
The factor of $\frac{1}{2}$ arises from the fact that the structures are non-plane, so that inside the brackets, each structure is obtained two times, one with its left and right subtrees of the root exchanged.

If there is a root gall, then the number of subtrees of the root gall is a value $k$ with $3 \leq k \leq n$. The total number of leaves in these subtrees is $n$, and the total number of galls is $g-1$ (the last gall is the root gall). As in enumerations for the unlabeled case in Section 4.2 of \cite{AgranatTamirEtAl24a}, because of the time-consistency condition, we have $k-2$ possible nodes in the root gall at which the reticulation node can be located. If the numbers of leaves in the subtrees are $(c_1, c_2, \ldots, c_k)$, then for each internal labeling of the subtrees, we have $\binom{n}{c_1, c_2, \ldots, c_k}$ options for distributing the $n$ labels among the subtrees. The numbers of leaves in the subtrees can be assigned according to each composition $(c_1, c_2, \ldots, c_k)$ in the set of compositions $C(n,k)$ of $n$ into $k$ parts; for each such composition, the numbers of galls in the subtrees are assigned by $(d_1-1,d_2-1,\ldots,d_k-1)$, where $(d_1,d_2,\ldots,d_k)$ is a composition of $g-1+k$ into $k$ parts. Hence, the contribution of the case with a root gall is 
\begin{equation}
\label{eq:engrec2}
     \frac{1}{2}\bigg{[}\sum_{k=3}^{n}(k-2)\sum_{\mathbf{c}\in C(n,k)}\sum_{\mathbf{d}\in C(g-1+k,k)}\binom{n}{c_1,c_2,\ldots, c_k}\prod_{i=1}^k e_{c_i,d_{i}-1}\bigg{]}.
\end{equation}
Again, the expression includes a factor of $\frac{1}{2}$ because the structures are nonplane, and each structure is obtained twice inside the brackets.

Summing (\ref{eq:engrec1}) and (\ref{eq:engrec2}), we have
\begin{align}
    \label{eq:engrec}
    e_{n,g} &= \frac{1}{2}\bigg{[}\bigg(\sum_{m=1}^{n-1} \binom{n}{m} \sum_{\ell=0}^g e_{m,\ell}e_{n-m,g-\ell}\bigg) \nonumber \\
    & + \bigg(\sum_{k=3}^{n}(k-2)\sum_{\mathbf{c}\in C(n,k)}\sum_{\mathbf{d}\in C(g-1+k,k)}\binom{n}{c_1,c_2,\ldots,c_k}\prod_{i=1}^k e_{c_i,d_{i}-1}\bigg)\bigg{]}.
\end{align}

Table \ref{table:gallnum} gives the numbers of labeled time-consistent galled trees calculated using the recursion for small $n$ and $g$. The total number of such trees across all $g$ accords with Theorem 8 of \cite{CardonaAndZhang20}. The case of $g=0$ recovers the familiar numbers of labeled trees with $n$ leaves. For the case of $g=1$, the numbers match those in \cite{Zhang19} who showed how to transition a tree-child network or a normal network with $n-1$ leaves and $k-1$ reticulations into a tree-child network or a normal network with $n$ leaves and $k$ reticulations. Specifically, \cite{Zhang19} showed that the number of labeled normal networks with one reticulation and $n$ labeled leaves---equivalent to labeled time-consistent galled trees with one gall---is ${(n+2)(2n)!}/(2^nn!)-3\cdot2^{n-1}n!$ for $n\geq3$.
\begin{table}[tb]
    \centering
    \small
    \begin{tabular}{|c|r|rrrrr|}
    \hline
Number of & \multicolumn{1}{c|}{Total number} & \multicolumn{5}{c|}{Number of trees with a fixed number of galls ($e_{n,g}$)} \\ 
leaves ($n$) & \multicolumn{1}{c|}{of trees ($a_n$)} & $g=0$ & $g=1$ & $g=2$ & $g=3$ & $g=4$ \\ \hline
 1 &              1 &          1 &             - &             - &              - &           - \\
 2 &              1 &          1 &             - &             - &              - &           - \\
 3 &              6 &          3 &             3 &             - &              - &           - \\
 4 &             69 &         15 &            54 &             - &              - &           - \\
 5 &          1,050 &        105 &           855 &            90 &              - &           - \\
 6 &         20,025 &        945 &        14,040 &          5040 &              - &           - \\
 7 &        464,310 &     10,395 &       248,535 &       197,820 &          7,560 &           - \\
 8 &     12,709,305 &    135,135 &     4,787,370 &     6,917,400 &        869,400 &           - \\
 9 &    401,112,810 &  2,027,025 &   100,361,835 &   233,859,150 &     63,617,400 &   1,247,400 \\
10 & 14,338,565,325 & 34,459,425 & 2,282,912,100 & 7,927,227,000 &  3,850,723,800 & 243,243,000 \\
\hline
\end{tabular}
\caption{Numbers of labeled time-consistent galled trees with specified numbers of leaves and galls. Entries $e_{n,g}$ are computed recursively (\ref{eq:engrec}).}
\label{table:gallnum}
\end{table}


\section{Discussion}

As phylogenetic networks have become increasingly central to mathematical phylogenetic studies, the enumerative combinatorics of network classes has been considered for many types of phylogenetic networks (Table \ref{table:survey}). We have focused here on the enumerative combinatorics of a highly restricted class of networks, the time-consistent galled trees. We have provided new derivations for unlabeled time-consistent galled trees and new results for labeled time-consistent galled trees, comparing the asymptotics of the latter with a fixed number of galls to corresponding asymptotics for other network classes.

For the unlabeled time-consistent galled trees, we have provided a derivation, using the symbolic method of analytic combinatorics, of the generating function that enumerates time-consistent galled trees with $n$ unlabeled leaves and any fixed number of galls $g$ (\ref{eq:eg1})-(\ref{eq:eg3b}); we had previously derived the generating function in \cite{AgranatTamirEtAl24b} by a recursive approach. The new derivation is relatively simple for small numbers of galls; the derivation for arbitrary $g$ proceeds through the bivariate generating function---which was derived in \cite{AgranatTamirEtAl24a} using recursion as well---and is somewhat more involved. 

Following the symbolic method, we derived generating functions enumerating time-consistent galled trees with $n$ labeled leaves and a fixed number of galls $g$ (\ref{eq:1termL})-(\ref{eq:2termL}). Relying on our asymptotic work on \emph{unlabeled} time-consistent galled trees with fixed $g$~\citep{AgranatTamirEtAl24b}, we have found asymptotic approximations for the number of labeled time-consistent galled trees with $n$ leaves and $g$ galls (\ref{eq:eng}).

This latter analysis produces a curious result. Many network classes with small numbers of reticulations all have the same asymptotic growth in the number of leaves~(\ref{eq:asymptotic}), a speed that in some classes can be explained combinatorially~\citep{FuchsEtAl24}. The number of time-consistent galled trees, however, grows more slowly in its subexponential term. Many of the inclusions in Table \ref{table:inclusion} are negligible in the sense that the smaller class of networks asymptotically has the same size as the larger class; for the time-consistent galled trees, however, the number of networks is asymptotically smaller than several other classes.

The results augment three earlier studies on unlabeled time-consistent galled trees~\citep{MathurAndRosenberg23, AgranatTamirEtAl24b, AgranatTamirEtAl24a}. For the unlabeled case, they provide simpler derivations. For the labeled case, they quickly produce new results and a comprarison with other structures. More generally, they contribute to furthering the enumerative combinatorics of phylogenetic networks.

\vskip .5cm
\noindent {\bf Acknowledgments.} We acknowledge grant support from National Science Foundation grant BCS-2116322 and from National Science and Technology Council grant NSTC-113-2115-M-004-004-MY3.
{\small \bibliography{refs} }
\end{document}